\documentclass[11pt]{article}
\usepackage{graphicx}
\usepackage{amssymb}
\usepackage{wasysym}
\def\L{{\bf L}}

\def\N{{\mathcal N}}
\def\D{{\mathcal D}}

\def\ve{\varepsilon}

\def\vp{\varphi}

\def\Hat{\widehat}
\def\Tilde{\widetilde}
\def\transv{\cap\!\!\!\!|~ }
\def\I{{\mathcal I}}

\def\C{{\mathcal C}}
\def\U{{\mathcal U}}
\def\S{{\mathcal S}}

\def\implies{\Longrightarrow}
\def\iff{\Longleftrightarrow}

\def\ds{\displaystyle}
\def\sqr#1#2{\vbox{\hrule height .#2pt
\hbox{\vrule width .#2pt height #1pt \kern #1pt
\vrule width .#2pt}\hrule height .#2pt }}
\def\square{\sqr74}
\def\endproof{\hphantom{MM}\hfill\llap{$\square$}\goodbreak}
\def\bega{\begin{array}}
\def\enda{\end{array}}
\def\begi{\begin{itemize}}
\def\endi{\end{itemize}}
\def\T{{\mathbb T}}
\def\R{{\mathbb R}}
\def\Z{{\mathbb Z}}

\def\ov{\overline}
\def\forall{\hbox{for all }~}
\def\v{\vskip 1em}
\def\vs{\vskip 2em}

\def\be{\begin{equation}}
\def\beq{\begin{equation}}
\def\bel{\begin{equation}\label}
\def\eeq{\end{equation}}

\begin{document}
\title{\bf Generic Regularity of Conservative Solutions 
to a Nonlinear  Wave Equation}
\vs

\author{Alberto Bressan$^{(*)}$ and  Geng Chen$^{(**)}$, \\    \\
{\small (*) Department of Mathematics, Penn State University,
University Park, Pa.~16802, U.S.A.}\\
{\small (**) School of Mathematics,
Georgia Institute of Technology,
Atlanta, Ga.~30332, U.S.A.}
\\
\,
\\
{\small e-mails:~ bressan@math.psu.edu~,~gchen73@math.gatech.edu}}

\maketitle

\begin{abstract} %\n 
The paper is concerned with  
conservative solutions to the nonlinear  %, second order 
wave equation
$u_{tt} - c(u)\big(c(u) u_x\big)_x~=~0$.    
For an open dense set of $\C^3$ 
initial data, we
prove that the solution is piecewise smooth in the $t$-$x$ plane, while
the gradient $u_x$ can blow up along  finitely 
many characteristic curves.  
The analysis is based on a variable transformation
introduced in \cite{BZ}, which reduces the equation to a semilinear 
system with smooth coefficients, followed by an application of
Thom's transversality theorem.
\end{abstract}
\v
\section{Introduction}
\setcounter{equation}{0}

Consider the  quasilinear
second order wave equation
\bel{1.1}
u_{tt} - c(u)\big(c(u) u_x\big)_x~=~0\,,
\qquad\qquad t\in [0,T], ~~x\in\R\,.\eeq
On the wave speed $c$ we assume
\begi
\item[{\bf (A)}] The map
 $c:\R\mapsto \R_+$ is
 smooth and  uniformly positive.   The quotient   $c'(u)/c(u)$ is uniformly 
 bounded. Moreover,
 the following generic condition is satisfied:
 \bel{morse}  c'(u)~=~0\qquad\implies\qquad c''(u)~\not=~ 0.\eeq 
 \endi
 Notice that, by (\ref{morse}), the derivative $c'(u)$ 
 vanishes only at isolated points.

The analysis in \cite{BZ, BCZ} shows that, for any initial data 
\bel{1.2}
u(0,x)~=~u_0(x)\,,\qquad
u_t(0,x)~=~u_1(x)\,,
\eeq
with 
$u_0\in H^1(\R)$, $u_1\in\L^2(\R)$, the  Cauchy problem admits a unique
conservative solution $u=u(t,x)$, H\"older continuous in the $t$-$x$ plane.
We recall that conservative solutions satisfy an additional conservation 
law for the energy, so that the total energy
$${\mathcal E}(t)~=~{1\over 2} \int [u_t^2 +  c^2(u) u_x^2]\, dx$$
coincides with a constant for a.e.~time $t$.
A detailed construction of a global semigroup of these solutions, 
including more singular
initial data, was carried out in \cite{HR}.

In the present paper we study the structure of these solutions.
Roughly speaking, we prove that, for generic  smooth 
initial data
$(u_0,u_1)$, the solution is piecewise smooth. Its gradient 
$u_x$ blows up along  finitely many smooth curves  in the $t$-$x$ plane.  
Our main result is
\v
{\bf Theorem 1.} 
{\it Let the function $u\mapsto c(u)$ satisfy the assumptions
{\bf (A)} and let  $T>0$ be given.
Then there exists an open dense set of initial data 
$$\D~\subset~ \Big(\C^3(\R)\cap H^1(\R)\Big) \times\Big(\C^2(\R)\cap\L^2(\R)\Big)$$
such that, for $(u_0,u_1)\in\D$, the conservative solution 
$u=u(t,x)$ of (\ref{1.1})-(\ref{1.2}) 
is twice continuously differentiable in the complement of finitely many characteristic curves 
$\gamma_i$, within the domain $[0,T]\times \R$.
}
\v

For the scalar conservation law in one space dimension,  a well known 
result by Schaeffer \cite{S} shows that
generic solutions are piecewise smooth, with finitely many shocks on any bounded domain 
in the $t$-$x$ plane.
A similar result was proved by Dafermos and Geng \cite{DG}, for a special 
$2\times 2$ Temple class 
system of conservation laws.  It remains an outstanding open problem 
to understand whether generic solutions to 
more general  $2\times 2$ systems (such as the p-system of isentropic gas dynamics)
remain piecewise smooth, with  finitely many shock curves.

The proof in \cite{S} relies on the Hopf-Lax representation formula, 
while the proof in \cite{DG}
is based on the analysis of solutions along characteristics.
In the present paper we take a quite different approach, based on the representation
of solutions in terms of a semilinear system introduced in \cite{BZ}.
In essence, the analysis in \cite{BZ} shows that, after a suitable change of variables,
the quantities
$$w~\doteq~2\arctan\bigl(u_t +c(u) u_x\bigr),
\qquad\qquad z ~\doteq~2\arctan\bigl(u_t -c(u) u_x\bigr),$$
satisfy a semilinear system of equations, w.r.t.~new 
independent variables $X$, $Y$. See (\ref{2.26})--(\ref{4.2}) 
in Section~2 for details.
Since this system has smooth coefficients, starting with smooth 
initial data
one obtains a globally defined
smooth solution.  To recover the singularities of the 
solution $u$ of (\ref{1.1})
in the original $t$-$x$ plane, it now suffices to study the level sets
\bel{wzp}  \{w(X,Y)=\pi\} \,,\qquad\qquad  \{z(X,Y)=\pi\} \,.  \eeq
Since $w$ and $z$ are smooth,  the generic structure of these
level sets
can be analyzed by  techniques of singularity theory \cite{Bloom, D, DD, G, T},
relying on Thom's transversality theorem.   
One should be aware that, while the map 
$(X,Y)\mapsto (t,x,u,w,z)$ is smooth, 
the inverse map $(t,x)\mapsto (X,Y)$ can have singularities.   
This variable transformation is indeed the source of 
singularities in the solution  $u=u(t,x)$   of (\ref{1.1}).

The present work was motivated by a research 
program aimed at the construction of a distance which 
renders Lipschitz continuous the semigroup of conservative solutions
of (\ref{1.1}).
Toward this goal, one needs
a dense set of piecewise smooth paths of solutions, 
whose weighted length can be controlled in time.
In the final section of this paper we thus consider
a 1-parameter family of initial data 
$\lambda\mapsto (u_0^\lambda, u_1^\lambda)$, with
$\lambda\in [0,1]$. 
We show that it can be uniformly approximated by a second path
of initial data
$\lambda\mapsto (\tilde u_0^\lambda, \tilde u_1^\lambda)$,
such that the corresponding solutions $\tilde u^\lambda =\tilde 
u^\lambda(t,x)$ of (\ref{1.1})
are piecewise smooth in the domain $[0,T]\times\R$,
for all except at most finitely values of $\lambda\in [0,1]$.
An application of this result to
the construction  of a Lipschitz metric will appear 
in the forthcoming paper~\cite{BC}.

The remainder of the paper is organized as follows.
In Section~2 we review the variable change introduced in \cite{BZ} and 
derive the semilinear system used in the construction 
of conservative solutions
to (\ref{1.1}).   In  Section~3 we construct families of smooth solutions 
to the semilinear system, depending on parameters.   
By a transversality argument, in Section~4 we show that
for almost all of these solutions  
the level sets (\ref{wzp}) satisfy a  number  of generic properties.  
After these preliminaries, the proof of Theorem~1
is completed in Section~5.
Finally, in Section~6 we prove a theorem on generic regularity 
for 1-parameter family of solutions.

For the nonlinear equation (\ref{1.1}), the formation of
singularities in finite time was first studied in \cite{GHZ}.
Based on the representations \cite{BZ, BH},
a detailed asymptotic description of structurally stable
singularities is given in \cite{BHY},  for 
conservative as well as  dissipative solutions.

We conjecture that the  regularity property stated in Theorem~1
should also hold for generic dissipative solutions of (\ref{1.1}). 
However, in the dissipative case
the corresponding semilinear system derived in \cite{BH}
contains discontinuous
terms,
and smooth initial data do not yield globally 
smooth solutions. For this reason, 
the techniques used in this paper can
no longer be applied.  We remark that, at the present time,
the uniqueness and continuous dependence of 
dissipative solutions to (\ref{1.1}) has not yet been proved,
for general initial data $(u_0, u_1)\in H^1(\R)\times \L^2(\R)$.

\v
\section{Review of the main equations}
\setcounter{equation}{0}

Consider the variables
\beq
\left\{
\begin{array}{rcl}
R & \doteq  &u_t+c(u)u_x\,, \\
S & \doteq  &u_t-c(u)u_x\,,
\end{array} \right.\label{2.1}
\eeq
so that
\beq
u_t~=~{R+S\over 2}\,,\qquad\qquad u_x~=~{R-S\over 2c}\,.\label{2.2}
\eeq
For a smooth solution of (\ref{1.1}), these variables  satisfy
\bel{2.3}
\left\{
\begin{array}{rcl}
R_t-cR_x &= & {c'\over 4c}(R^2-S^2), \\ [3mm]
S_t+cS_x &= & {c'\over 4c}(S^2-R^2).
\end{array} \right.
\eeq
In addition,  $R^2$ and $S^2$ satisfy the  balance laws
\bel{2.4}\left\{
\begin{array}{rcl}
(R^2)_t - (cR^2)_x & = & {c'\over 2c}(R^2S - RS^2)\, , \\ [3mm]
(S^2)_t + (cS^2)_x & = &  {c'\over 2c}(S^2R-SR^2 )\,.
\end{array}
\right.
\eeq
As a consequence, for smooth solutions the following quantity is  conserved:
\beq
E~\doteq ~{1\over 2}\big(u_t^2+c^2u_x^2\big)~=~{R^2+S^2\over 4}
%\,,\qquad\qquad M~\doteq~ -u_tu_x~=~{S^2-R^2\over 4c}
\,.\label{2.5}
\eeq
One can think of $R^2$ and $S^2$ as the energy densities of backward and 
forward moving waves, respectively.  Notice that these are not separately conserved. 
Indeed, by (\ref{2.4}) energy can be exchanged between forward and backward waves.

It is well known that, even for smooth
initial data, the quantities $u_t, u_x$ can blow up in finite time
\cite{GHZ}.
To deal with possibly unbounded values of $R,S$, following \cite{BZ}
it is convenient to introduce a new set of dependent variables:
\bel{wzdef}
w~\doteq ~2\arctan R\,,\qquad\qquad z~\doteq~ 2\arctan S\,.
\eeq
Using (\ref{2.3}), we obtain the equations
\beq
w_t-c\,w_x~=~{2\over 1+R^2}(R_t-c\,R_x)
~=~{c'\over 2c}{R^2-S^2\over 1+R^2}\,, \label{2.10}
\eeq
\beq
z_t+c\,z_x~=~{2\over 1+S^2}(S_t+c\,S_x)
~=~{c'\over 2c}{S^2-R^2\over 1+S^2}\,. \label{2.11}
\eeq

\begin{figure}[htb]
\centering
\includegraphics[width=0.55\textwidth]{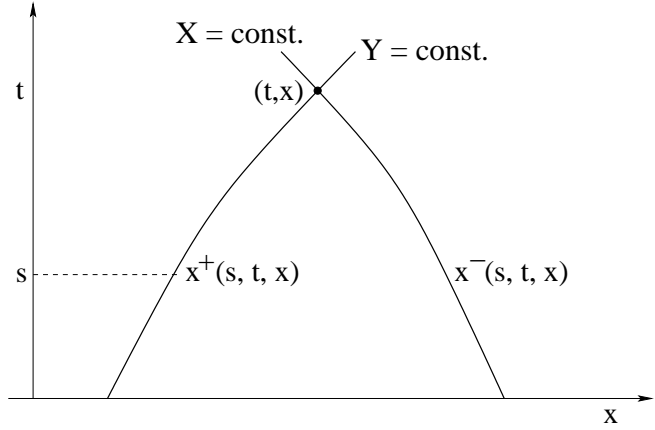} %AB's
\caption{\small Characteristic curves.  As new coordinates of the point $(t,x)$
we choose the values $(X,Y) = \bigl(x^-(0, t,x)\,,~-x^+(0,t,x)\bigr)$.}
\label{f:b77}
\end{figure}
We now
perform a further change of independent variables (Fig.~\ref{f:b77}).
Consider
the equations for the  backward and forward characteristics:
\beq
\dot x^-~=~-c(u)\,,\qquad\qquad \dot x^+~=~c(u)\,,\label{2.12}
\eeq
where the upper dot denotes a derivative w.r.t.~time. 
The characteristics passing through the point $(t,x)$
will be denoted by
$$
s~\mapsto~ x^-(s,t,x)\,,\qquad\qquad s~\mapsto ~x^+(s,t,x)\,,
$$
respectively.
As coordinates $(X,Y)$ of a point $(t,x)$ we shall use the intersections
of these characteristics with the $x$-axis, namely
\beq
X~\doteq~ x^-(0,t,x)
\,,\qquad\qquad Y~\doteq ~- x^+(0,t,x)\,.\label{XY}
\eeq
Of course this implies
\beq
X_t- c(u)X_x~ =~0\,,\qquad\qquad Y_t+ c(u)Y_x ~=~0\,,\label{2.14} 
\eeq
\beq
(X_x)_t- (c\,X_x)_x~=~0\,,\qquad\qquad (Y_x)_t+(c\,Y_x)_x~=~0\,.\label{2.15}
\eeq
For any smooth function $f$, using (\ref{2.14}) one finds
\beq\left\{
\begin{array}{ccccccr}
f_t+cf_x &=& f_XX_t+f_Y Y_t+cf_X X_x+cf_Y Y_x & = & (X_t+cX_x)f_X
& = & 2cX_x f_X\,,\cr\cr
f_t-cf_x &=& f_XX_t+f_Y Y_t-cf_X X_x-cf_Y Y_x & = & (Y_t-cY_x)f_Y
& = & -2cY_x f_Y\,.
\end{array}\right. \label{2.16} 
\eeq
We now introduce the further variables
\bel{2.17}
p~\doteq~ {1+R^2\over X_x}\,,\qquad\qquad q~\doteq ~{1+S^2\over -Y_x}\,.
\eeq
Notice that the above definitions imply
\beq
{1\over X_x}~=~ {p\over 1+R^2}~=~p\,\cos^2{w\over 2}\,,\qquad\qquad
{-1\over Y_x}~=~{q\over 1+S^2}~=~q\,\cos^2{z\over 2}\,.\label{2.18}
\eeq

Starting with the nonlinear equation (\ref{1.1}),
using $X,Y$ as independent variables one obtains a semilinear
hyperbolic system with smooth coefficients for the variables $u,w,z,p,q$, namely
\beq
\left\{
\begin{array}{ccc}
u_X &= & {\sin w\over 4c} \, p\,,\\[4mm]
 u_Y &= & {\sin z\over 4c} \, q\,,
\end{array}\right. \label{2.26} %{1.18}
\eeq
\beq
\left\{
\begin{array}{ccc}
w_Y&=&{c'\over 8c^2}\,( \cos z - \cos w)\,q\,,\\[4mm]
z_X&=&{c'\over 8c^2}\,( \cos w - \cos z)\,p\,,
\end{array}\right. \label{2.24} %{1.16}
\eeq
\beq
\left\{
\begin{array}{ccc}
p_Y &= & {c'\over 8c^2}\, \big(\sin z-\sin w\big)\,pq\,,\\[4mm]
q_X &= & {c'\over 8c^2}\, \big(\sin w-\sin z\big)\,pq\,.
\end{array}\right.\label{2.25} %{1.17}
\eeq
The map $(X,Y)\mapsto (t,x)$ can  be constructed as follows.
Setting $f=x$, then $f=t$  in the two equations at
(\ref{2.16}), we find
$$
\left\{
\begin{array}{rcr}
c &=& 2cX_x\,x_X\,, \\[3mm]
-c&=& -2cY_x\,x_Y\,,
\end{array}\right.\qquad\qquad
\left\{
\begin{array}{rcr}
1&= &2cX_x\,t_X\,,\\[3mm]
1&=& -2cY_x\,t_Y\,,
\end{array}\right.
$$
respectively.
Therefore, using (\ref{2.18}) we obtain
\beq
\left\{
\begin{array}{ccccr}
x_X&~=~&{1\over 2X_x}&~=~&{(1+\cos w)\,p\over 4}\,,\\[4mm]
x_Y&~=~&{1\over 2Y_x}&~=~&-{(1+\cos z)\,q\over 4}\,,
\end{array}\right. \label{4.1}
\eeq
\beq
\left\{
\begin{array}{ccccr}
 t_X&=&{1\over 2cX_x}&=&{(1+\cos w)\,p\over 4c}\,,\\[4mm]
t_Y&=&{1\over -2cY_x}&=&{(1+\cos z)\,q\over 4c}\,.
\end{array}\right. \label{4.2}\eeq
\v
Given the initial data (\ref{1.2}), the corresponding  boundary 
data for   (\ref{2.24})-(\ref{4.2}) can be determined as follows.
In the $X$-$Y$ plane, consider the line 
$$\gamma_0~=~\{X+Y=0\}~\subset~\R^2$$
parameterized as
%\bel{gpar}
$x~\mapsto ~(\ov X(x), \,\ov Y(x))~\doteq~(x,\, -x)$.
Along $\gamma_0$ we can assign the boundary data
$(\ov u,\ov w, \ov z, \ov p, \ov q)$ by setting
 \beq{\ov u }~ =~u_0(x)\,,\qquad\qquad
\left\{
\begin{array}{rcl}
{\ov w } &= & 2\arctan R(0,x)\,,\\ {\ov z} &= & 2\arctan S(0,x)\,,
\end{array}\right.\qquad \qquad
\left\{
\begin{array}{rcl}
{\ov p} &\equiv & 1+ R^2(0,x)\,,\\
{\ov q} &\equiv & 1+S^2(0,x)\,,
\end{array}\right.  \label{2.28}
\eeq
at each point $(x,-x)\in\gamma_0$.
We recall that, at time $t=0$, by (\ref{1.2}) one has
$$\bega{lr}R(0,x) &=~(u_t + c(u) u_x)(0,x) ~=~
u_1(x) + c(u_0(x)) u_{0,x}(x),\\[4mm]
S(0,x) &=~(u_t - c(u) u_x)(0,x) ~=~u_1(x) - c(u_0(x)) u_{0,x}(x).\enda$$
\v
{\bf Remark 1.}
Since the semilinear system (\ref{2.24})--(\ref{4.2}) has smooth coefficients,
for smooth initial data all components of the solution remain smooth on the entire $X$-$Y$
plane.   As proved in \cite{BZ}, 
the quadratic terms in (\ref{2.25}) (containing the product
$pq$)  account
for transversal wave interactions and 
do not produce finite time blow up of the
variables $p,q$.   Moreover, if the values of $p,q$ are uniformly positive along a line $\{X+Y=\kappa\}$, then 
they remain uniformly positive
on compact sets of the $X$-$Y$ plane.  Throughout this paper, we always
consider solutions of (\ref{2.24})--(\ref{4.2}) where $p,q>0$.
\v
By 
expressing the solution $u(X,Y)$ in terms of the original variables
$(t,x)$, one obtains a solution of the Cauchy problem
(\ref{1.1})-(\ref{1.2}). Indeed, the following was proved in \cite{BZ}.
\v
{\bf Lemma 1.} {\it Let $(u,w,z,p,q,x,t)$ be a smooth solution
to the system  (\ref{2.26})--(\ref{4.2}), with $p,q>0$. 
Then the set of points
\bel{graph}\Big\{ \bigl(t(X,Y), \, x(X,Y), 
\, u(X,Y)\bigr)\,;~~(X,Y)\in\R^2\Big\}\eeq
is the graph of a conservative solution to 
the variational wave equation (\ref{1.1}).}
\v

We observe that, while the functions
\bel{Ldef}
(X,Y)~\mapsto~ u(X,Y),\qquad\qquad (X,Y)~\mapsto ~\Lambda(X,Y) ~\doteq
~\bigl(t(X,Y), \, x(X,Y)
\bigr)\eeq
are globally smooth,
the map $\Lambda:\R^2\mapsto\R^2$ may not have a smooth inverse.  
Indeed, $\Lambda$ may not even 
be one-to-one.  Therefore, the solution $u(t,x) = u(\Lambda^{-1}(t,x))$
can fail to be smooth.  This happens precisely at points where
the Jacobian matrix $D\Lambda$ is not invertible.   
By (\ref{4.1})-(\ref{4.2}), singularities occur when  
$\cos w= -1$  or  $\cos z= -1$.
\v
{\bf Remark 2.} The system (\ref{2.24})--(\ref{4.2}) is
invariant under translation by $2\pi$
in $w$ and $z$.  We can thus think of 
$w,z$ as points in the  quotient manifold 
$\T=\R/2\pi \Z$. 
Throughout the following we take advantage of this fact and 
regard a solution of (\ref{2.24})--(\ref{4.2}) as a 
map $(X,Y)\mapsto(u,w,z,p,q,x,t)$ from
$\R^2$ into $\R\times\T\times\T \times\R\times\R\times\R$.
Observe that we have the implications
\bel{coswz}\bega{l}
w~\not=~\pi\qquad\implies\qquad \cos w ~>~ -1\,,\\[4mm]
z~\not=~\pi\qquad\implies\qquad \cos z ~>~ -1\,.\enda\eeq
\v
{\bf Remark 3.}  In general,  many distinct solutions to the system
(\ref{2.26})--(\ref{4.2}) can yield the same solution $u=u(t,x)$ of (\ref{1.1}).

Indeed, let $(u,w,z,p,q,x,t)(X,Y)$ be one particular solution. 
Let $\phi,\psi:\R\mapsto\R$ be two $\C^2$ bijections, 
with $\phi'>0$ and $\psi'>0$.  Introduce the new independent and dependent variables
$(\Tilde X,\Tilde Y)$ and $(\tilde u,\tilde w,\tilde z,\tilde p,\tilde q,\tilde x,\tilde t)$
by setting
\bel{TXY}
X~=~\phi(\Tilde X)\,,\qquad\qquad Y~=~\psi(\Tilde Y),\eeq
\bel{TUWZ}
(\tilde u,\tilde w,\tilde z,\tilde x, \tilde t)(\Tilde X,\Tilde Y)~=~(u,w,z,p,q,x,t)(X,Y),\eeq
\bel{TPQ}\left\{\bega{rl}
\tilde p(\Tilde X,\Tilde Y)&=~p(X,Y)\cdot \phi'(\Tilde X),\\[4mm]
\tilde q(\Tilde X,\Tilde Y)&=~q(X,Y)\cdot \psi'(\Tilde Y).\enda\right.\eeq
Then, as functions of $(\Tilde X, \Tilde Y)$,  the
variables $(\tilde u,\tilde w,\tilde z,\tilde p,\tilde q,\tilde x,\tilde t)$
provide another solution of the same system (\ref{2.26})--(\ref{4.2}). 
Moreover, by (\ref{TUWZ}) the set
\bel{graph2}\Big\{ \Bigl(\tilde t(\Tilde X,\Tilde Y), \, \tilde x(\Tilde X,\Tilde Y), 
\,\tilde u(\Tilde X,\Tilde Y)\Bigr)\,;~~
(\Tilde X,\Tilde Y)\in\R^2\Big\}\eeq
coincides with the set (\ref{graph}). Hence it is the graph of the same solution 
$u$ of (\ref{1.1}).  
One can regard the variable transformation (\ref{TXY})
simply as a relabeling of forward and backward characteristics, in the solution $u$. A detailed analysis of relabeling symmetries, in 
connection with the Camassa-Hom equation, can be found in \cite{GHR}.

For future reference we observe that
$$\tilde w_{\Tilde X}(\Tilde X,\Tilde Y)~=~w_X(X,Y)\cdot \phi'(\Tilde X)\,,$$
$$\tilde w_{\Tilde X\Tilde X}(\Tilde X,\Tilde Y)~=~w_{XX}(X,Y)\cdot [\phi'(\Tilde X)]^2
+w_X(X,Y) \cdot \phi''(\Tilde X)\,.$$
In particular, one has the equivalences
\bel{degen}
\bega{rl}
\tilde w_{\Tilde X}(\Tilde X,\Tilde Y)~=~0\qquad &\iff \qquad 
w_X(X,Y)~=~0,\\[4mm]
\tilde z_{\Tilde Y}(\Tilde X,\Tilde Y)~=~0\qquad &\iff \qquad 
z_Y(X,Y)~=~0,\\[4mm]
(\tilde w_{\Tilde X},\,\tilde w_{\Tilde X\Tilde X})(\Tilde X,\Tilde Y)~=~(0,0)\qquad 
&\iff \qquad 
(w_X, w_{XX})(X,Y)~=~(0,0),\\[4mm]
(\tilde z_{\Tilde Y},\,\tilde z_{\Tilde Y\Tilde Y})(\Tilde X,\Tilde Y)~=~(0,0)\qquad 
&\iff \qquad 
(z_Y, z_{YY})(X,Y)~=~(0,0)\,.\enda\eeq

\v
\subsection{Compatible boundary data}
More generally,
instead of (\ref{2.28}) we can assign  boundary data 
for the system (\ref{2.16})--(\ref{2.25}) on a  line $\gamma=\{ X+Y=\kappa\}$.
Namely:
 \bel{bdata}  u(s,\kappa-s)~ =~\ov u(s)\,,\qquad\qquad
\left\{
\begin{array}{rl}
w(s,\kappa-s)&= ~\ov w(s)\,,\\ z(s,\kappa-s)&= ~ \ov z(s)\,,
\end{array}\right.\qquad \qquad
\left\{
\begin{array}{rcl}
p(s,\kappa-s)&=~ \ov p(s)\,,\\
q(s,\kappa-s)&=~\ov q(s)\,,
\end{array}\right.  
\eeq
for suitable smooth functions $\ov u, \ov w,\ov z,\ov p,\ov q$.  
If both identities in (\ref{2.26}) hold,
then 
\bel{cc}
{d\over ds}  \ov u(s) ~=~ {d\over ds } u(s,\kappa-s)~=~u_X - u_Y~=~
{\sin w\over 4c} \,p- {\sin z\over 4c}\,q\,.\eeq
The boundary data should thus satisfy 
the compatibility condition
\bel{cc1} {d\over ds}  \ov u(s)~=~{\sin \ov w(s)\over 4c(\ov u(s))} \,
\ov p(s)- {\sin \ov z(s)\over 4c(\ov u(s))}\,\ov q(s)\eeq

As remarked earlier,  the system (\ref{2.26})--(\ref{2.25}) is overdetermined.
Indeed, the function $u=u(X,Y)$ could be recovered by either one of the 
identities in (\ref{2.26}).    We now prove that, if the compatibility condition
(\ref{cc1}) holds, then any smooth solution
satisfying one of the identities in (\ref{2.26}) satisfies the other as well.
\v
{\bf Lemma 2.} {\it 
Let $u,w,z,p,q$ be smooth functions on $\R^2$ which satisfy 
(\ref{2.24})-(\ref{2.25})
together with  the boundary conditions (\ref{bdata}) 
along the line $\gamma=\{ X+Y=\kappa\}$.  Assume that the compatibility 
condition (\ref{cc1}) holds. 
Then one has 
\bel{uY}u_Y~=~{\sin z\over 4c(u)}\,q  \qquad\qquad 
\forall (X,Y)\in \R^2\eeq
if and only if
\bel{uX}u_X~=~{\sin w\over 4c(u)}\,p \qquad\qquad 
\forall (X,Y)\in \R^2.\eeq
}
\v
{\bf Proof.} Consider the smooth, strictly increasing function
$$\Phi(u)~=~\int_0^u 4c(s)\, ds\,.$$
 Observe that the identities  (\ref{uY}),  (\ref{uX}) are equivalent
 respectively to
\bel{PXY}\Phi(u)_Y~=~\sin z\cdot q\,,\qquad\qquad
\Phi(u)_X~=~\sin w\cdot p\,. \eeq

 Assume  that (\ref{uY}) holds.  Then
\bel{P5}\Phi(u( X, Y))~=~\Phi(u(X, \,\kappa- X))+
 \int_{\kappa- X}^{ Y}  
[\sin z\cdot q](X,  s)\,  ds\,.\eeq
Differentiating w.r.t.~$X$,  and using the first equations
in (\ref{2.24})-(\ref{2.25}) together with the
compatibility condition (\ref{cc1}) we obtain
\bel{PX}
\bega{rl}\Phi(u)_X(X,Y)
&\ds=~\Phi'(u)\cdot[u_X-u_Y](X, \kappa-X)+ 
[\sin z\cdot q](X,  \kappa-X)
\cr\cr
&\ds\qquad\qquad+ \int_{\kappa-Y}^Y [ \cos z\cdot z_X q + 
\sin z\cdot q_X](X,s)\, ds
\cr\cr
&\ds=~[\sin w\cdot p](
X, \kappa-X) +\int_{\kappa-Y}^Y 
\Big[{c'(u)\over 8 c^2(u)}\bigl(1+\sin(z+w)\bigr)pq\Big]
(X,s)\, ds\cr\cr
&\ds=~[\sin w\cdot p](
X, \kappa-X) +\int_{\kappa-Y}^Y {\partial\over\partial Y}
[\sin w\cdot p]
(X,s)\, ds~=~[\sin w\cdot p](X,Y)\,.
\enda
\eeq
We have thus proved that the second identity in (\ref{PXY}) holds.
This is equivalent to (\ref{uX}).

The converse implication is proved in the same way.
\endproof

%\begin{figure}[htb]
%\centering
%\includegraphics[width=0.55\textwidth]{FIG/wa69.eps} %AB's
%\caption{\small .}
%\label{f:wa69}
%\end{figure}
\v
Next, consider initial data for $t,x$, 
on the curve $\gamma=\{X+Y=\kappa\}$, say
\bel{bd2}  x(s,\kappa-s)~ =~\ov x(s)\,,\qquad\qquad
t(s,\kappa-s)~=~\ov t(s)\,. 
\eeq
Using (\ref{4.1})-(\ref{4.2}) we derive 
the  compatibility conditions
\bel{cc5}{d\over ds} \ov x(s) ~=~{(1+\cos \ov q(s)) \ov p(s) 
+ (1+\cos \ov z(s) )\ov q(s)\over 4}\,,\eeq
\bel{cc6}
{d\over ds} \ov t(s) ~=~{(1+\cos \ov  w(s)) \ov p(s) 
- (1+\cos \ov z(s) )\ov q(s)
\over 4 c(\ov u(s))}\,.\eeq
\v
{\bf Lemma 3.} {\it 
Let $(u,w,z,p,q)(X,Y)$ be a solution of the system  
(\ref{2.26})--(\ref{2.25}).  Then there exists a solution 
$(t,x)(X,Y)$ of (\ref{4.1})-(\ref{4.2}) with boundary 
data (\ref{bd2}) if and only if the compatibility
conditions 
(\ref{cc5})-(\ref{cc6}) are satisfied.}
\v
{\bf Proof.} {\bf 1.} Assume that 
 the equations (\ref{4.1})-(\ref{4.2}) are satisfied
 for all $(X,Y)\in\R^2$. In particular, they are satisfied
 along the curve $\gamma=\{X+Y=\kappa\}$.
 This implies
$${d\over ds}\, \ov x(s)~=~
{d\over ds}\, \ov x(s, \kappa-s)~=~[x_X-x_Y](s,\kappa-s)~
=~{(1+\cos  w)  p 
+ (1+\cos  z ) q\over 4}\,,$$
where the right hand side is evaluated at $(X,Y) = (x, \kappa-s)$.
Hence (\ref{cc5}) holds.   The identity (\ref{cc6}) 
is derived in the same way.
\v
{\bf 2.} Next, assume that the compatibility conditions 
(\ref{cc5})-(\ref{cc6}) are satisfied.   To prove that 
 (\ref{4.1}) admits a solution,
it then suffices to check 
that the differential form
$${(1+\cos w)p\over 4}\, dX -{(1+\cos z)q\over 4}\, dY$$
is closed.
This is true because
\bel{xXY}
\begin{array}{rl}\ds
\left[ {(1+\cos w)p\over 4}\right]_Y
&= ~\ds-{\sin w\over 4} w_Y\,p+{1+\cos w\over 4}\,p_Y\\[4mm]
&\ds=~ -{\sin w\over 4} \cdot {c'\over 8 c^2}(\cos z - \cos w)\, pq + {1+\cos w
\over 4} \cdot {c'\over 8c^2}(\sin z-\sin w) \, pq\\[4mm]
&\ds=~  {c'\over 32 c^2}\Big[ 
(1+\cos w)\sin z-(1+\cos z)\sin w\Big]\, pq~=~
 \left[ -{(1+\cos z)q\over 4}
\right]_X\,.
\end{array} 
\eeq
Similarly, to prove that (\ref{4.2}) admits a solution,
it suffices to check that the 
differential form 
$${(1+\cos w)p\over 4c}\, dX +{(1+\cos z)q\over 4c}\, dY$$
is closed.
This is true because 
\bel{tXY}
\begin{array}{l}\ds
\left[ {(1+\cos w)p\over 4c}\right]_Y
\ds= ~-{\sin w\over 4c} w_Y\,p - {1+\cos w\over 4c^2} c' u_Y\,p
+{1+\cos w\over 4c}\,p_Y\\[4mm]
=~\ds -{\sin w\over 4c} \cdot {c'\over 8 c^2}(\cos z - \cos w)\, pq 
- {1+\cos w\over 4c^2} c' {\sin z\over 4c}\,pq
+ {1+\cos w
\over 4c} \cdot {c'\over 8c^2}(\sin z-\sin w) \, pq\\[4mm]
=~\ds - {c'\over 32 c^3}\Big[ (1+\cos w)\sin z+(1+\cos z)\sin w\Big]\, pq~=~
 \left[ {(1+\cos z)q\over 4c}\right]_X \,.
\end{array} 
\eeq
\endproof
\v
\v
{\bf Remark 4.}   Let a solution $(u,w,z,p,q)$ of (\ref{2.26})--(\ref{2.25})
be given.
If we assign the values of $t,x$ 
at a single point $(X_0, Y_0)$, then by the compatibility conditions (\ref{cc5})-(\ref{cc6})
and the equations (\ref{4.1})-(\ref{4.2}) the functions $t(X,Y)$, $x(X,Y)$ are uniquely determined for all $(X,Y)\in\R^2$.  Choosing different 
values of $t,x$  at the point $(X_0, Y_0)$  we obtain the same solution
 $u=u(t,x)$ of  (\ref{1.1}), up to a shift of coordinates in the  $t,x$ plane.  
 \v
\section{Families of perturbed solutions}
\setcounter{equation}{0}
 Let a point $(X_0, Y_0)$ be given and consider the
line 
\bel{gd}\gamma~\doteq~\Big\{ (X,Y)\,;~~~X+Y=\kappa\Big\},
\qquad\qquad \kappa \doteq X_0+Y_0\,.\eeq
We can then arbitrarily assign the values of $w,z,p,q$
at every point $(X,Y)\in \gamma$.   Moreover, we can arbitrarily
choose the values of $u,x,t$ at the single point $(X_0, Y_0)$.
In turn, these choices uniquely determine functions
$u,x,t$ on $\gamma$ which satisfy the compatibility 
conditions (\ref{cc1}) and (\ref{cc5})-(\ref{cc6}).

Based on this observation, we can construct several families
of perturbations
of a given solution of (\ref{2.26})--(\ref{2.25}).
The main goal of this section is to prove
\v
{\bf Lemma 4.} {\it Let the assumption {\bf (A)} hold.
%Consider initial data
%$(u_0, u_1)\in H^1(\R)\cap\L^2(\R)$ in (\ref{1.2}), 
%with $u_0,u_1\in\C^\infty$.
Let $(u, w,z, p, q)$ be a smooth solution of the system
(\ref{2.26})-(\ref{2.25}) and  
let a point $(X_0, Y_0)\in\R^2$ be given.
 
 \begi
 \item[{\bf (1)}] If $(w,w_X,w_{XX})(X_0,Y_0)=(\pi,0,0)$, 
 then there exists a 3-parameter family of smooth 
 solutions $(u^\theta, w^\theta,z^\theta,  p^\theta, q^\theta)$ 
 of (\ref{2.26})-(\ref{2.25}),
  depending smoothly on $\theta\in\R^3$,
 such that the following holds.
\begi
\item[(i)] When $\theta=0\in\R^3$ one recovers the original solution,
namely $(u^0, w^0,z^0, p^0, q^0)=(u, w,z, p, q)$.
 
\item[(ii)] At the point
 $(X_0,Y_0)$, when $\theta=0$ one has
 \bel{pert1}\hbox{\rm rank }D_\theta (w^\theta\,,~
w^\theta_X\,,~w^\theta_{XX})~=~3\,.\eeq
\endi

\item[{\bf (2)}] If $(w,z,w_X)(X_0,Y_0)=(\pi,\pi,0)$, then 
there exists a 3-parameter family of smooth 
 solutions $(u^\theta, w^\theta,z^\theta,  p^\theta, q^\theta)$
  satisfying (i)-(ii) as above, with (\ref{pert1})
 replaced by
 \bel{pert2}\hbox{\rm rank }D_\theta (w^\theta\,,~
z^\theta\,,~w^\theta_X)~=~3\,.\eeq

\item[{\bf (3)}] If $(w,w_X, c'(u))(X_0,Y_0)=(\pi,0,0)$,  
then there exists a 3-parameter family of smooth 
 solutions $(u^\theta, w^\theta,z^\theta,  p^\theta, q^\theta)$ satisfying (i)-(ii) as above, with (\ref{pert1})
 replaced by
\bel{pert3}\hbox{\rm rank }D_\theta \bigl(w^\theta\,,~
w^\theta_X\,,~c'(u^\theta)\bigr)~=~3\,.\eeq

\endi
}
\v
For example (\ref{pert1}) means that we can construct  perturbed solutions, depending on
parameters $\theta_1,\theta_2,\theta_3$,  such that the Jacobian matrix
\bel{mat0}D_\theta (w^\theta\,,~
w^\theta_X\,,~w^\theta_{XX})~=~\left(\bega{ccc} {\partial\over\partial\theta_1} w 
& {\partial\over\partial\theta_2} w
& {\partial\over\partial\theta_3} w\cr\cr {\partial\over\partial\theta_1} w_X& {\partial\over\partial\theta_2} w_X
& {\partial\over\partial\theta_3} w_X\cr\cr
{\partial\over\partial\theta_1} w_{XX}& {\partial\over\partial\theta_2} w_{XX}
& {\partial\over\partial\theta_3} w_{XX}\enda\right),\eeq
computed at $\theta=0$, has full rank at the point $(X_0, Y_0)$.
\v
\subsection{ Proof of Lemma 4}
Let $(u,w,z,p,q)$ be a $\C^\infty$ solution of the system
(\ref{2.26})-(\ref{2.25}).   
Given the point $(X_0, Y_0)$, 
consider the line $\gamma$  in (\ref{gd}) and 
let $(\ov u,\ov w,\ov z,\ov p,\ov q)$
be the values of the solution along $\gamma$, as in 
(\ref{bdata}).

For future use, we compute the values of $w_X, w_{XX}$
at the point $(X_0, Y_0)$.   At any point $(s, \kappa-s)\in 
\gamma$ we have
$$w_X- w_Y~=~\ov w'(s)\,,
\qquad\qquad z_X- z_Y~=~\ov z'(s)\,,\qquad\qquad
q_X- q_Y~=~\ov q'(s)\,.$$
Here and in the sequel, a prime denotes derivative w.r.t.~the 
parameter $s$ along the curve $\gamma$.
Using  (\ref{2.24})-(\ref{2.25})
 we obtain
\bel{wX}
w_X(X_0, Y_0)~ =~\ov w' + {c'(\ov u)\over 8c^2(\ov u)}
(\cos \ov z-\cos \ov w)\ov q\,,\eeq
\bel{zY}z_Y(X_0, Y_0)~=~-\ov z'+ {c'(\ov u)\over 8c^2(\ov u)}
(\cos \ov w - \cos \ov z) \ov p\,,\eeq
\bel{qY}q_Y(X_0, Y_0)~=~-\ov q'+ {c'(\ov u)\over 8c^2(\ov u)}
(\sin \ov w - \sin \ov z) \ov p\ov q\,,\eeq
where all terms on the right hand sides are evaluated at $s=X_0$.

A further differentiation yields
$${d^2\over ds^2} \ov w(s)~=~{d\over ds} 
\bigl[w_X(s, \kappa-s)-w_Y(s,\kappa-s)\bigr]~=~[w_{XX}+w_{YY}-2w_{XY}]
(s,\kappa-s).$$
Using (\ref{2.26})--(\ref{2.25}) together with (\ref{wX})--(\ref{qY})
we obtain
\bel{wXY}
\bega{l}
w_{YX}(X_0, Y_0)\cr\cr
\quad =~
\left({c'(\ov u)\over 8c^2(\ov u)}\right)' {\ov u}_X\,( \cos\ov z - \cos\ov w)\,\ov q
+
{c'\over 8c^2}\,({\ov w}_X \sin\ov w- {\ov z}_X \sin \ov z)\,\ov q
+
{c'\over 8c^2}\,( \cos \ov z - \cos \ov w)\,{\ov q}_X\cr\cr
\quad =~\left({c'\over 8c^2}\right)' \frac{\sin \ov w}{4c}\,( \cos\ov z - \cos\ov w)\,\ov p \ov q
\cr\cr\qquad\qquad
+{c'\over 8c^2}\,\left\{\big(\ov w' + {c'\over 8c^2} (\cos \ov z-\cos \ov w)\ov q\big) \sin\ov w- {c'\over 8c^2}
(\cos \ov w - \cos \ov z) \ov p\sin \ov z\right\}\,\ov q
\cr\cr\qquad\qquad 
+({c'\over 8c^2})^2\,( \cos \ov z - \cos \ov w)\,
(\sin \ov w - \sin \ov z) \ov p\ov q
\cr\cr \quad \doteq ~f_1
\enda
\eeq
and
\bel{wYY}
\bega{l}
w_{YY}(X_0, Y_0)\cr\cr
\quad =~
\left({c'(\ov u)\over 8c^2(\ov u)}\right)'\ov u_Y\,( \cos\ov z - \cos\ov w)\,\ov q
+
{c'\over 8c^2}\,(\ov w_Y \sin \ov w- \ov z_Y \sin\ov  z)\,\ov q
+
{c'\over 8c^2}\,( \cos\ov z - \cos\ov w)\,\ov q_Y\cr\cr
\quad =~\left({c'\over 8c^2}\right)' \frac{\sin \ov z}{4c}\,( \cos\ov z - \cos\ov w)\,\ov q^2
\cr\cr
\qquad\qquad
+{c'\over 8c^2}\,\left\{ {c'\over 8c^2} (\cos \ov z-\cos \ov w)\ov q\sin\ov w- \big(-\ov z'+{c'\over 8c^2}
(\cos \ov w - \cos \ov z) \ov p\big)\sin \ov z\right\}\,\ov q
\cr\cr
\qquad\qquad +
{c'\over 8c^2}\,( \cos \ov z - \cos \ov w)\,
\big(-\ov q'+ {c'\over 8c^2}
(\sin \ov w - \sin \ov z) \ov p\ov q\big)\cr\cr
\quad \doteq~ f_2\,.
\enda
\eeq
Hence
\bel{wXX}
w_{XX}(X_0, Y_0)~=~\ov w'' + 2f_1 - f_2\,.\eeq
\v
We  now  construct families $(\ov u^\theta,
\ov q^\theta,\ov z^\theta, \ov p^\theta, \ov q^\theta)$
of perturbations of the data
(\ref{bdata}) along the curve $\gamma$, so that 
at the point $(X_0, Y_0)$ the matrices in (\ref{pert1})--(\ref{pert3})
have full rank.  
These perturbations will have the form
\bel{p00}\left\{\bega{rl}
\ov w^\theta(s) &=~\ov w(s) + \sum_{i=1}^3\theta_i W_i(s)\,,\cr\cr
\ov z^\theta(s) &=~\ov z(s) + \sum_{i=1}^3\theta_i Z_i(s)\,,\enda
\right.
\qquad\qquad 
\left\{\bega{rl}
\ov p^\theta(s) &=~\ov p(s) + \sum_{i=1}^3\theta_i P_i(s)\,,\cr\cr
\ov q^\theta(s) &=~\ov q(s) + \sum_{i=1}^3\theta_i Q_i(s)\,,\enda
\right.\eeq
for suitable functions $W_i, Z_i, P_i, Q_i\in\C^\infty_c$.
Moreover, at the point $s=X_0$ we set
\bel{u00}
\ov u^\theta(X_0) ~=~\ov u(X_0) + \sum_{i=1,2,3}\theta_i U_i\,.\eeq
In turn, the above definitions together
with the compatibility conditions (\ref{cc1}) determine
 the values of $\ov u^\theta(s)$ for all $s\in\R$.
In particular, for each $\theta\in\R^3$ we obtain 
a unique solution of
the semilinear system (\ref{2.26})--(\ref{2.25}).

We observe that the functions $W_i, Z_i, P_i,Q_i$ can be chosen arbitrarily. Hence
at the point $s=X_0$, $\theta=0$, 
we can arbitrarily assign  all derivatives
$${d\over d\theta} {d^k\over ds^k}\ov w^\theta,
\qquad {d\over d\theta}{d^k\over ds^k} \ov z^\theta,
\qquad
{d\over d\theta}{d^k\over ds^k} \ov p^\theta\,,
\qquad{d\over d\theta} {d^k\over ds^k} \ov q^\theta\,,$$
with $k=0,1,2,\ldots$~   Moreover, we can
arbitrarily choose the quantity
${d\over d\theta} \ov u^\theta(X_0)$,
while  all higher order derivatives ${d\over d\theta}{d^k\over ds^k} 
\ov u^\theta$,
with $k\geq 1$, are then determined by the 
compatibility condition (\ref{cc1}).
\v
{\bf 1.} 
To achieve (\ref{pert1}), we choose
perturbations $(\ov u^{\theta_i},
\ov w^{\theta_i},\ov z^{\theta_i}, \ov p^{\theta_i}, 
\ov q^{\theta_i})$, $i=1,2,3$, so that
the Jacobian matrix of first order derivatives
w.r.t.~$\theta_1,\theta_2,\theta_3$, computed  
at $s= X_0$ and $\theta=0$, is given by
\[D_\theta \left(\bega{c}
\ov u\cr
\ov w\cr
\ov z\cr
\ov z'\cr
\ov w'\cr
\ov w''\cr
\ov p\cr
\ov q\cr
\ov q'\enda\right)~=~
\left(\bega{ccc}0&0&0\cr
1&0&0\cr
0&0&0\cr
0&0&0\cr
0&1&0\cr
0&0&1\cr
0&0&0\cr
0&0&0\cr
0&0&0
\enda\right)\,.\]
At the point $(X_0, Y_0)$, by (\ref{wX}) and (\ref{wXX})  this yields
\[D_\theta \left(\bega{c}
 w\cr  w_X\cr w_{XX}\enda\right)~=~
\left(\bega{ccc}
1&0&0\cr
*&1&0\cr
*&*&1\cr
\enda\right)\,.\]
Notice that, for the third family of perturbations (corresponding to the third column),  the first order variations of 
$f_1$ and $f_2$ in (\ref{wXY})-(\ref{wYY}) 
both vanish at $(X_0, Y_0)$.
This achieves (\ref{pert1}).
\v
{\bf 2.} To achieve (\ref{pert2}), we choose
perturbations $(\ov u^{\theta_i},
\ov w^{\theta_i},\ov z^{\theta_i}, \ov p^{\theta_i}, 
\ov q^{\theta_i})$, $i=1,2,3$, so that  at $s= X_0$ and $\theta=0$
one has
\bel{DT1}D_\theta \left(\bega{c}
\ov u\cr
\ov w\cr
\ov z\cr\ov w'\cr
\ov p\cr\ov q\enda\right)~=~
\left(\bega{ccc}0&0&0\cr
1&0&0\cr
0&1&0\cr
0&0&1\cr
0&0&0\cr
0&0&0
\enda\right)\,.\eeq
At the point $(X_0, Y_0)$, by (\ref{wX}) this yields
\bel{DT2}D_\theta \left(\bega{c}
 w\cr  z\cr w_X\enda\right)~=~
\left(\bega{ccc}
1&0&0\cr
0&1&0\cr
*&*&1\cr
\enda\right)\,.\eeq
Hence (\ref{pert2}) holds.
\v
{\bf 3.} Finally, we construct three families of 
perturbations satisfying (\ref{pert3}). If at $(X_0, Y_0)$
we have $c'(u(X_0,Y_0)) = 0 $, then  the assumption
{\bf (A)} implies
\bel{c''}
%K~\doteq~\left({c'(u)\over 8c^2(u)}\right)'~=~{c''(u)\over 8c^2(u)}~ \neq~ 0.
c''\bigl(u(X_0, Y_0)\bigr)~\not= ~0\,.
\eeq
To achieve (\ref{pert3}), we choose  three families of 
perturbations such that, at $s= X_0$ and $\theta=0$, 
\bel{DT5}D_\theta \left(\bega{c}
\ov u\cr
\ov w\cr
\ov z\cr\ov w'\cr
\ov p\cr\ov q\enda\right)~=~
\left(\bega{ccc}0&0&1\cr
1&0&0\cr
0&0&0\cr
0&1&0\cr
0&0&0\cr
0&0&0
\enda\right)\,.\eeq
At the point $(X_0, Y_0)$, by (\ref{wX}) and the first equation in (\ref{2.24}), this yields
\bel{DT6}D_\theta \left(\bega{c}
 w\cr  w_X\cr c'(u)\enda\right)~=~
\left(\bega{ccc}
1&0&0\cr
*&1&*\cr
*&0&c''(u)\cr
\enda\right)\,\eeq
This achieves (\ref{pert3}). 
\endproof

\v
\section{Generic solutions of the semilinear system}
\setcounter{equation}{0}

In this section we study smooth solutions to the semilinear 
system (\ref{2.24})--(\ref{4.2}),  determining the generic structure
of the level sets $\{(X,Y)\,;~~w(X,Y) =\pi\}$ and $\{(X,Y)\,;~~z(X,Y)=\pi\}$. 
%In turn, this will provide a proof of Theorem~1.
\v
{\bf Lemma 5.} {\it Let the function $u\mapsto c(u)$ satisfy the assumptions
{\bf (A)} and consider a compact domain  of the form
\bel{GTdef}\Gamma~\doteq~\Big\{(X,Y)\,;~~|X|+|Y|~\leq~M\Big\}.\eeq  
Call $\S$ the family of all 
 $\C^2$ solutions to
the system (\ref{2.26})--(\ref{2.25}), with $p,q>0$ for all 
$(X,Y)\in\R^2$. Moreover, call $\S'\subset\S$ the subfamily
of all solutions $(u,w,z,p,q)$ 
 such that, for $(X,Y)\in \Gamma$, none of the following 
values is attained:
\bel{never1}\left\{\bega{rl} (w,w_X,w_{XX}) &=~(\pi,0,0),\\[3mm]
(z,z_Y,z_{YY}) &=~(\pi,0,0),\enda\right.\eeq
\bel{never2}\left\{\bega{rl} (w,z,w_X) &=~(\pi,\pi,0),\\[3mm]
(w,z,z_Y) &=~(\pi,\pi,0),\enda\right. \eeq
\bel{never3}
\left\{\bega{rl} (w,w_X,c'(u)) &=~(\pi,0,0),\\[3mm]
(z,z_Y, c'(u)) &=~(\pi,0,0).\enda\right.\eeq
Then $\S'$ is a relatively open and dense subset of $\S$, in the topology induced by $\C^2(\Gamma)$.
}

\begin{figure}[htbp]
\centering
\includegraphics[width=0.7\textwidth]{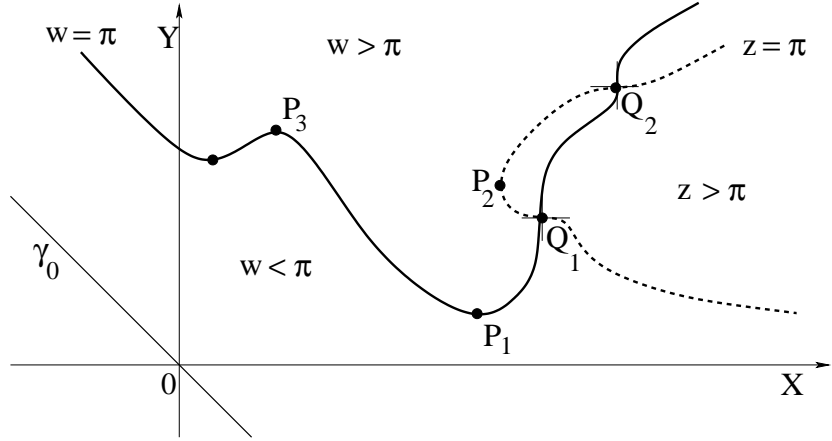}
\caption{ \small Two level sets $\{w=\pi\}$ and $\{z=\pi\}$,
in a generic solution of (\ref{2.26})--(\ref{2.25}). 
At $P_1$, $P_2$ one has $w=\pi$, $w_X=0$ 
while the generic conditions imply $w_Y\not=0$, 
$w_{XX}\not= 0$.    At the points
 $Q_1$, $Q_2$  where the two singular curves cross, by (\ref{2.24}) one has $w_Y=z_X=0$, while the generic conditions imply
 $w_X\not= 0$, $z_Y\not= 0$. Hence the two curves have a perpendicular intersection. }
\label{f:wa67}
\end{figure}

Some words of explanation are in order (Fig.~\ref{f:wa67}). 
Asking that the values in (\ref{never1}) are never attained 
is equivalent to the implications
\bel{impl1}\left\{\bega{cl}
w~=~\pi\quad\hbox{and}\quad w_X~=~0\qquad &\implies\qquad w_{XX}~\not= 0,\\[3mm]
z~=~\pi\quad\hbox{and}\quad z_Y~=~0\qquad &\implies\qquad z_{YY}~\not= 0.\enda
\right.
\eeq
Writing the level curves in the form $\{w(X,Y)=\pi\}= \{ Y=\vp(X)\}$ and 
$\{z(X,Y)=\pi\}= \{ X=\psi(Y)\}$, this imposes some restrictions at the points
where $\vp'=0$ or $\psi'=0$.
%This imposes some restrictions at the points where the level curves
%$\{w=\pi\}$ and $\{z=\pi\}$ initiate or terminate.

Asking that the values in (\ref{never2}) are never attained 
is equivalent to the implication
\bel{impl2}
[w~=~\pi\quad\hbox{and}\quad z~=~\pi]\qquad \implies\qquad 
[w_X~\not= ~0\quad\hbox{and}\quad z_Y~\not=~0].
\eeq
This imposes restrictions at points where two level curves
$\{w=\pi\}$ and $\{z=\pi\}$ cross each other.
%More generally, by the same arguments one can
%construct a perturbed solution $(u,w,z,p,q)$ for which the values
%\bel{never7}
%\left\{\bega{cl} (w,\, \cos z,\,w_X) &=~(\pi,-1,0),\\[3mm]
%(z,\,\cos w,\,z_Y) &=~(\pi,-1,0),\enda\right.\eeq
%are never attained, for any $(X,Y)\in\Gamma$.

Finally, the lemma states the existence of a perturbed solution
such that values (\ref{never3}) 
are never attained.  To understand the meaning of this condition,
consider a solution  which never attains any of the values 
in (\ref{never2})-(\ref{never3}). In this case, by 
(\ref{coswz}) the 
conditions
$w=\pi$ and $w_X=0$ together imply
$$w_Y~=~{c'(u)\over 8c^2(u)}(\cos z+1) q~\not= ~0\,.$$
This is equivalent to  the implication 
$$w~=~\pi~\qquad\implies\qquad (w_X, w_Y)~\not=~(0,0).$$
By the implicit function theorem, 
the level set $\{w=\pi\}$ is then the union of regular curves in the $X$-$Y$
plane (restricted to the domain $\Gamma$). 
Similarly,  the level set $\{z=\pi\}$ will be a 
union of regular curves.

\begin{figure}[htbp]
\centering
\includegraphics[width=0.6\textwidth]{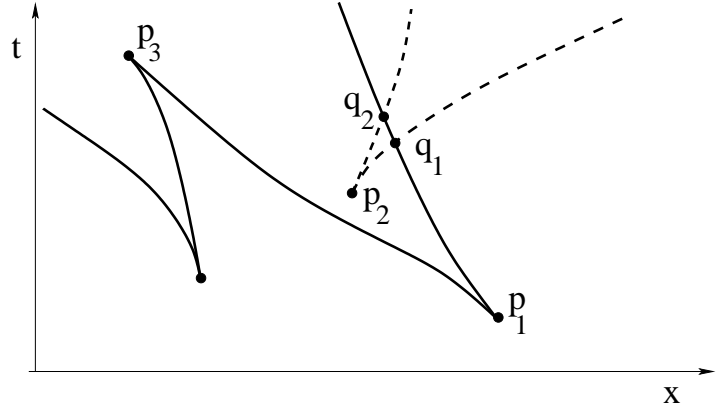}
\caption{ \small   A generic solution $u=u(t,x)$ 
of (\ref{1.1}) with smooth initial data
remains smooth outside finitely many singular points and finitely many singular curves where
$u_x\to\pm\infty$.  The curves where $u$ is singular are the images of the curves where
$w=\pi$ or $z=\pi$ in Fig.~\ref{f:wa67}, under the map $(t,x)=\Lambda(X,Y)$
at (\ref{Ldef}).
Here $p_i=\Lambda(P_i)$ are  points  where singular curves originate or terminate, while
$q_j=\Lambda(Q_j)$
are points where two singular curves cross.}
\label{f:wa68}
\end{figure}

We shall  give a proof of Lemma 5,
using Lemma~4 together with Thom's transversality theorem. 
For readers' convenience, we first review some basic definitions \cite{Bloom, GG, T}.
\v
{\bf Definition (map transverse to a submanifold).} 
Let $f : X \mapsto  Y$ be a smooth map of manifolds and let 
$W$ be  a submanifold of $Y$. We say that 
$f$ is {\it transverse} to $W$ at  a point $p\in X$, 
and write $f \transv_p W$, if
\begi
\item either $ f ( p ) \notin  W$,
\item or else $f(p) \in  W$
 and $T_{f(p)}Y = (df)_p(T_pX) + T_{f(p)}W$.
 \endi
 Here $T_p X$ denotes the tangent space to $X$ at the point $p\in X$, 
while $T_qY$ and $T_q W$ denote respectively the tangent spaces to $Y$ and to $W$ at the point  
 $q\in W\subset Y$.  Finally, $(df)_p:T_pX\mapsto T_{f(p)}Y$ denotes the differential
 of the map $f$ at the point $p$.

 We say that $ f$  is {\it transverse} to $W$, and write $f \transv W$, if $f \transv_p W$ for every $p\in X$.

In the special case where $W=\{y\}$ consists of a single point, $f \transv W$ 
if and only if $y$ is a  regular value of $f$, in the following sense. 

\v
{\bf Definition (regular value).}   Let $f : X \mapsto  Y$ be a smooth map of manifolds.
A point $y\in Y$ is a {\it regular value} if,
for every $p\in X$ such that $f(p)=y$, one has
$$T_{f(p)}Y~ =~ (df)_p(T_pX).$$

\v
{\bf Transversality Theorem.} {\it Let $X$, $\Theta$, and $Y$ be smooth manifolds, 
$W$ a submanifold
of $Y$. 
Let $\theta\mapsto \phi^\theta$ be a smooth map which to each 
$\theta\in \Theta$ associates a
function  $ \phi^\theta\in  \C^\infty(X,Y)$, and define 
$\Phi:X\times \Theta\mapsto Y$  by setting 
$\Phi(x,\theta)=\phi^\theta(x)$. 
If $\Phi\transv W$ then the set 
$\{\theta\in \Theta,;~~\phi^\theta\transv W\}$ 
is dense in $\Theta$.}

For a proof, see \cite{Bloom, GG}.
\v
\subsection{Proof of Lemma 5.}

{\bf 1.} We shall use the representation
\bel{Srep} 
S'~=~S_1\cap S_2\cap S_3\cap S_4\cap S_5\cap S_6\,,\eeq
where $S_1,\ldots,S_6\subset\S$ are the families of   solutions
for which one of  the six values listed in 
(\ref{never1})--(\ref{never3}) is never attained on $\Gamma$. 
For example, 
$\S_1$ is the set of all solutions such that 
\bel{S1}(w,w_X,w_{XX})(X,Y) ~\not=~(\pi,0,0)\qquad\qquad 
\forall (X,Y)\in\Gamma,\eeq
while 
$\S_6$ is the set of all solutions such that 
\bel{S6}(z, z_Y, c'(u))(X,Y) ~\not=~(\pi,0,0)\qquad\qquad 
\forall (X,Y)\in\Gamma.\eeq

Since $\Gamma$ is a compact domain, 
it is clear that each $\S_i$ is a relatively
open subset of $\S$, in the topology of $\C^2(\Gamma)$.
In the remainder of the proof we will show that each $S_i$ is dense
on $\S$.
\v
{\bf 2.}  Let $(u,w,z,p,q)$ be any $\C^2$ solution of 
(\ref{2.26})--(\ref{2.25}), with $p,q>0$.  By a smooth approximation
of the  data along the line $\gamma_0=\{X+Y=0\}$, it is not 
restrictive to assume that $u,w,z,p,q\in\C^\infty(\R^2)$.
We begin by looking at the first condition in (\ref{never1}).

Given any  point $(X_0, Y_0)\in \Gamma$, two cases can occur.

CASE 1: $(w, w_X, w_{XX})(X_0, Y_0)~\not=~(\pi, 0,0)$.
In this case, by continuity, there exists a neighborhood $\N$ of $(X_0, Y_0)$
in the $X$-$Y$ plane where we still have $(w, w_X, w_{XX})~\not=~(\pi, 0,0)$.

CASE 2: $(w, w_X, w_{XX})(X_0, Y_0)~=~(\pi, 0,0)$.  In this case,
by Lemma~4 we can find a 3-parameter family of 
solutions $(u^\theta, w^\theta, z^\theta, p^\theta, q^\theta)$
such that the $3\times 3$ Jacobian matrix of the map 
\bel{5map}
(\theta_1,\theta_2,\theta_3)~\mapsto~\Big(w^\theta(X,Y)\,, ~w_X^\theta(X,Y)\,,w_{XX}^\theta(X,Y)\Big)\eeq
has rank 3 at the point $(X_0, Y_0)$, when $\theta=0$.  
By continuity, this matrix still has rank 3
on a neighborhood $\N$ of $(X_0,Y_0)$, for $\theta$ sufficiently
close to zero.

We now choose finitely many points $(X_i, Y_i)$, $i=1,\ldots,n$, 
such that the corresponding open neighborhoods
$\N_{(X_i,Y_i)}$ cover the compact set $\Gamma$. 
Call $n_\I$ the cardinality of the set of indices
\bel{Idef}
\I~\doteq~\bigl\{ i\,;~~(w, w_X, w_{XX})(X_i, Y_i)~=~(\pi,0,0)\bigr\}
\eeq
so that CASE 2 applies, and set $N=3n_\I$.
\v
{\bf 3.}
Let $\Omega\supset\Gamma$ be an open set contained in the union
of the  neighborhoods $\N_{(X_i,Y_i)}$, and 
call $B_\ve\doteq\bigl\{ \theta\in \R^N\,;~~|\theta|<\ve\bigr\}$
 the open ball of radius $\ve$  in $\R^N$.  

We shall construct a family $(u^\theta, w^\theta,z^\theta, 
p^\theta, q^\theta)$ of smooth solutions to (\ref{2.26})--(\ref{2.25}),
 such that the map
\bel{Nmap}
(X,Y,\theta)~\mapsto~\Big(w^\theta(X,Y)\,, ~w_X^\theta(X,Y)\,,~w_{XX}^\theta(X,Y)\Big)\eeq
from $\Omega\times B_\ve$ 
into $\R^3$ has $(\pi,0,0)$ as a regular value.
Toward this goal, we need
to combine  perturbations based at possibly different points
$(X_i, Y_i)$ into a single $N$-parameter family
of perturbed solutions.

Let $(u,w,z,p,q)(X,Y)$ be a solution to the system  
(\ref{2.26})--(\ref{2.25}).
For each $k=1,\ldots,N$,
let a point $(X_k,Y_k)$ be given,
together with  a number  $U_k\in \R$ and functions
$W_k, Z_k, P_k, Q_k\in \C^\infty_c(\R)$.
By the previous analysis,  a 1-parameter family of perturbed solutions
to (\ref{2.26})--(\ref{2.25}) is then determined as follows.
For $|\ve|<\ve_k$ sufficiently small,  let 
\bel{solt}
(u^\ve, w^\ve, z^\ve, p^\ve, q^\ve) ~\doteq~\Psi_k^\ve(u,w,z,p,q)\eeq
be the unique solution of (\ref{2.26})--(\ref{2.25}) with 
data assigned on the
line $\gamma_k\doteq \{X+Y=X_k+Y_k\}$  by setting
$$u^\ve(X_k, Y_k) ~=~u(X_k, Y_k)+\ve U_k\,,$$
while for $(X,Y)\in \gamma_k$
$$w^\ve ~= ~w+\ve W_k\,,\qquad 
z^\ve ~= ~z+\ve Z_k\,,\qquad 
p^\ve ~= ~p+\ve P_k\,,\qquad 
q^\ve ~= ~q+\ve Q_k\,.$$
Given $(\theta_1,\ldots,\theta_{N})$, a perturbation of the 
original solution $(u,w,z,p,q)$ is defined as the composition
of $N$ perturbations:
\bel{pert5}
(u^\theta,w^\theta,z^\theta,p^\theta,q^\theta)~\doteq~
\Psi_{N}^{\theta_{N}}\circ\cdots\circ\Psi_1^{\theta_1}(u,w,z,p,q).\eeq
\v
{\bf 4.}
At each point $(X_i, Y_i)$ with $i\in\I$, we can apply Lemma~4 
 and obtain
three 1-parameter families of perturbed solutions 
so that the Jacobian matrix (\ref{5map}) has rank 3 on 
$\N_{(X_i,Y_i)}$,
for all $\theta$ small enough.

Combining all these perturbations,
we  obtain an $N$-parameter family 
of solutions  such that 
the value 
$(\pi,0,0)$ is a regular value
for  the map  (\ref{Nmap}), from $\Omega\times B_\ve$ into $\R^3$.

By the transversality theorem, for a.e.~$\theta$ the value 
$(\pi,0,0)$ is a regular value
for  the map $(X,Y)~\mapsto~
\Big(w^\theta(X,Y)\,, ~w_X^\theta(X,Y)\,,~w_{XX}^\theta(X,Y)\Big)$
 from $\Omega$ into $\R^3$.    
Since $\Omega$ has dimension 2, for a.e.~$\theta$ the corresponding solution
 $(u^\theta, w^\theta, z^\theta, p^\theta, q^\theta)$ 
 has the property that
$$\bigl(w^\theta(X,Y)\,, ~w_X^\theta(X,Y)\,,~w_{XX}^\theta(X,Y)
\bigr)~\not=~(\pi,0,0)$$
for all $(X,Y)\in \Gamma$.
This proves that the set $\S_1$ of solutions for which (\ref{S1})
holds is dense on $\S$.
\v
{\bf 5.} Repeating the above construction, we obtain that 
each $\S_i$, $i=1,\ldots,6$,  is a relatively open, dense 
subset of $\S$.  By (\ref{Srep}), the intersection $S'$ is 
is a relatively open, dense 
subset of $\S$.
\endproof

\section{Proof of Theorem 1.}
\setcounter{equation}{0}
Consider the product space
\bel{Udef}\U~\doteq~
\Big(\C^3(\R)\cap H^1(\R)\Big) 
\times\Big(\C^2(\R)\cap\L^2(\R)\Big)\eeq
with norm
$$\bigl\| (u_0,u_1)\bigr\|_\U~\doteq~\|u_0\|_{\C^3} + \|u_0\|_{H^1} + \|u_1\|_{\C^2}
+ \|u_1\|_{\L^2}\,.$$
Given  initial data $(\hat u_0, \hat u_1)~\in~ \U$,
consider the open ball 
\bel{Bd}B_\delta~\doteq~\Big\{ 
(u_0, u_1)\in\U\,;~~
\bigl\|(u_0, u_1)-(\hat u_0,\hat u_1)\bigr\|_\U~<~\delta\Big\}.
\eeq
Theorem 1 will be proved by showing that, for any 
$(\hat u_0, \hat u_1)~\in~ \U$ there exists a radius $\delta>0$
and an open dense subset 
$\Hat\D\subseteq B_\delta$, with the following property:
For every initial data
$(u_0,u_1)\in\Hat \D$, the conservative solution 
$u=u(t,x)$ of (\ref{1.1})-(\ref{1.2}) 
is twice continuously differentiable in the 
complement of finitely many characteristic curves 
$\gamma_i$, within the domain $[0,T]\times \R$.
\v
{\bf 1.}   Let $(\hat u_0, \hat u_1)~\in~ \U$ be given.
By the definition of the space $\U$ in (\ref{Udef}), 
as $|x|\to\infty$ we have
\bel{00}
\hat u_0(x)\to 0,\qquad \hat u_{0,x}(x)\to 0,
\qquad \hat u_1(x)\to 0.\eeq
Hence  the corresponding functions $R,S$ in (\ref{2.1})
satsfy
$$R(0,x)~\to ~0,\qquad\qquad S(0,x)~\to~0.$$
{}From (\ref{2.3}), it follows that the functions $R,S$ remain 
uniformly bounded on a domain of the form
$\{ (t,x)\,;~~t\in [0,T],~~~|x|\geq r\}$, for $r$ sufficiently 
large.  More generally, we can choose
$\delta>0$ such that,
for every initial data $(u_0, u_1)\in B_\delta$,
the corresponding solution $u(t,x)$ remains 
twice continuously differentiable
 on the outer domain
\bel{outer}\Big\{ (t,x)\,;~~t\in [0,T],~~~|x|\geq \rho\Big\},\eeq
for some $\rho>0$ sufficiently large.  Its singularities
can thus occur only on the compact domain 
$[0,T]\times [-\rho, \,\rho]$.

The subset $\Hat D\subset B_\delta$ is now defined as follows.
~$(u_0, u_1)\in \Hat D$ if $(u_0, u_1)\in B_\delta$
and moreover, for the corresponding solution
$(u,w,z,p,q)$ of (\ref{2.26})--(\ref{4.2}) with boundary data  
(\ref{2.28}), the values (\ref{never1})--(\ref{never3}) are never attained, 
for any $(X,Y)$
such that 
\bel{bad}(t(X,Y), \, x(X,Y))~\in~[0,T]\times [-\rho,\,\rho].\eeq
It is important to observe that, by (\ref{degen}), 
the above condition is independent
of the relabeling (\ref{TXY}).  
\v
{\bf 2.}   For any $(u_0, u_1)\in B_\delta$ we now consider the
corresponding solution $(t,x,u,w,z,p,q)$ of the system 
(\ref{2.26})--(\ref{4.2}), with boundary data as in (\ref{2.28}).
Let $\Lambda$ be the map  at (\ref{Ldef}) and let
$\Gamma$ be the square with side $2M$ in the $X$-$Y$ plane, 
as in (\ref{GTdef}).

By choosing  $M$ large enough, and by possibly shrinking the 
radius $\delta$, we can achieve the inclusion
\bel{sur}[0,T]\times [-\rho,\rho]~\subset~
\Lambda(\Gamma)\,,\eeq
for every $(u_0, u_1)\in B_\delta$.

\v
{\bf 3.}  We begin by proving that $\Hat \D$ is open, in the 
topology of $\C^3\times\C^2$.
Indeed, consider initial data $(u_0,u_1)\in\Hat \D$ and let 
$(u_0^\nu, u_1^\nu)_{\nu\geq 1}$
be a sequence of initial data converging to $(u_0, u_1)$.

Assume, by contradiction, that $(u_0^\nu, u_1^\nu)\notin \Hat \D$ for all $\nu\geq 1$.
To fix the ideas, let $(X^\nu, Y^\nu)$ be points at which 
the corresponding solutions $(u^\nu, w^\nu, z^\nu, p^\nu, q^\nu)$ satisfy
\bel{bad1}
(w^\nu, w^\nu_X, w^\nu_{XX}) (X^\nu, Y^\nu)~=~(\pi,\, 0,\,0),\qquad\quad (t^\nu,x^\nu)
(X^\nu, Y^\nu)~\in~[0,T]\times [-\rho,\,\rho],\eeq
for all $\nu\geq 1$.    
By (\ref{sur}), 
since the domain $\Gamma$ in (\ref{GTdef}) is compact,
by possibly taking a subsequence we can assume 
$(X^{\nu}, Y^{\nu})\to (\ov X, \ov Y)$.   By continuity, this
implies 
$$(w, w_X, w_{XX})(\ov X, \ov Y) ~=~(\pi, \,0,\,0),
\qquad\qquad (t,x)
(\ov X,\, \ov Y)~\in~[0,T]\times [-\rho,\,\rho],$$
contradicting the  assumption $(u_0, u_1)\in \Hat\D$.

The other cases in (\ref{never1})--(\ref{never3}) are handled in the same way.
This proves that $\Hat\D$ is open.
\v
{\bf 4.}  Next, we claim that $\Hat D$ is dense in $B_\delta$.
Indeed, let $(u_0,u_1)\in B_\delta$ be given.
By an arbitrarily small perturbation (measured in the norm of $\U$),
we can assume that $u_0, u_1\in\C^\infty$.

Using Lemma~5, we can construct a sequence of solutions 
$(u^\nu, w^\nu, z^\nu, p^\nu, q^\nu, x^\nu, t^\nu)$ of (\ref{2.26})--(\ref{4.2})
such that: 
\begi
\item[(i)] For every bounded set $\Omega\subset\R^2$ and any $k\geq 1$,
the $\C^k$ norm of the difference satisfies
\bel{lim0}
\lim_{\nu\to\infty}~\Big\|(u^\nu-u, ~w^\nu-w, ~z^\nu-z, ~p^\nu-p, 
~q^\nu-q, ~x^\nu-x, ~t^\nu-t)
\Big\|_{\C^k(\Omega)}~=~0.\eeq 
\item[(ii)] For every $\nu\geq 1$, 
the values in (\ref{never1})--(\ref{never3}) are never attained, 
for any $(X,Y)\in \Gamma$. 
\endi
Consider the   corresponding solutions 
$u^\nu(t,x)$ of (\ref{1.1}), with graph 
$$\Big\{ \bigl(u^\nu(X,Y),~ t^\nu(X,Y),~ x^\nu(X,Y)\bigr)\,;~~(X,Y)\in\R^2\Big\}
~\subset~\R^3.$$
For $t=0$, by (\ref{lim0}) the corresponding sequence of initial values satisfies
\bel{l5}\lim_{\nu\to\infty}~\bigl\| u^\nu(0,\cdot) - u_0\bigr\|_{\C^k([a,b])}~=~0,
\qquad\qquad 
\lim_{\nu\to\infty}~\bigl\| u^\nu_t(0,\cdot) - u_1\bigr\|_{\C^k([a,b])}~=~0,\eeq
for every bounded interval $[a,b]$.

Next, consider a cutoff function $\eta\in \C^\infty_c$ such that
\bel{eta}
\bega{l}\eta(x)~=~1\qquad\hbox{if}~~|x|\leq r,\\[4mm]
\eta(x)~=~0\qquad\hbox{if}~~|x|\geq r+1,\enda\eeq
with  $r>\!>\rho$ sufficiently large. 
For every $\nu\geq 1$, consider the initial data
$$\tilde u^\nu_0~\doteq~\eta u_0^\nu+(1-\eta)u_0\,,
\qquad  \tilde u^\nu_1~\doteq~\eta u_1^\nu+(1-\eta)u_1\,.$$
By (\ref{l5}) we have
\bel{l6}\lim_{\nu\to\infty}~\bigl\| (\tilde u_0^\nu - u_0, 
~\tilde u_1^\nu-u_1)\bigr\|_\U~=~0.\eeq
Moreover, if $r>0$ was chosen large enough, we have
$$\tilde u^\nu(t,x)~=~u^\nu(t,x)\qquad\qquad\forall 
(t,x)\in [0,T]\times [-\rho,\rho]\,,$$
while $\tilde u^\nu$ remains $\C^2$ on the outer domain (\ref{outer}).
%\bel{outer} \bigl\{(t,x)\,;~~t\in [0,T]\,,~~|x|>\rho\bigr\}.\eeq
The above implies $(\tilde u^\nu_0, \tilde u^\nu_1)\in \Hat \D$
for all $\nu\geq 1$ sufficiently large,
proving that $\Hat \D$ is dense
on $B_\delta$.
\v
{\bf 5.} To complete the proof we need to show that,
for every initial data $(u_0, u_1)\in \Hat\D$, the 
solution $u(t,x)$ of (\ref{1.1}) is piecewise $\C^2$ on the 
domain  $[0,T]\times\R$.

By the previous arguments, we already know that $u$ 
is $\C^2$ on the outer domain (\ref{outer}).  It thus remains
to study the singularities of $u$ on the inner domain 
$[0,T]\times [-\rho,\,\rho]$.
For this purpose,
call $(u,w,z,p,q,t,x)(X,Y)$ the corresponding solution of 
(\ref{2.26})--(\ref{4.2}), with boundary data as in (\ref{2.28}). 
By (\ref{sur}), every point of the inner domain is contained in the
image of the square $\Gamma$ in (\ref{GTdef}).
 
Consider a point $(X_0,Y_0)\in\Gamma$.   Two cases can occur.

CASE 1: $w(X_0,Y_0)\not= \pi$ and $z(X_0,Y_0)\not= \pi$.   
By (\ref{4.1})-(\ref{4.2}) it follows
$$\det\left(\bega{cc}x_X&x_Y\cr t_X&t_Y\enda\right)~=~
{(1+\cos w)p\over 4}\cdot {(1+\cos z)p\over 4c}+{(1+\cos z)q\over 4}\cdot {(1+\cos w)p\over 4c}~>~0.$$
Hence the map $(X,Y)\mapsto (x,t)$ is locally invertible 
in a neighborhood of $(X_0,Y_0)$.  We can thus conclude that 
the function $u$ is $\C^2$ in a neighborhood of the point 
$\bigl(t(X_0, Y_0),
x(X_0,Y_0)\bigr)$.
 
CASE 2: $w(X_0, Y_0)=\pi$.  In this case we have either $w_X(X_0,Y_0)\not= 0$, or else by (\ref{2.24})
\bel{wY0}w_Y(X_0,Y_0)~=~{c'(u)\over 8 c^2(u)}(\cos z +1)q~\not=~0.\eeq
Indeed, we always have  $c(u)>0$ and $q>0$. 
Moreover, by construction the values $(w,z,w_X)= (\pi,\pi,0)$
and $(w,w_X, c'(u))=(\pi,0,0)$ are never attained in $\Gamma$. 
This implies (\ref{wY0}).

By the implicit function theorem, we thus conclude that the sets
\bel{SWZ}S^w~\doteq~\{(X,Y)\in \Gamma\,;~~w(X,Y)=\pi\},\qquad
\qquad S^z~\doteq~\{(X,Y)\in \Gamma\,;~~z(X,Y)=\pi\}\eeq
are the union of finitely many $\C^2$ curves.

The set of points $(t,x)$ where $u$ is singular 
coincides with the image of
the two sets $S^w$, $S^z$ under the $\C^2$ map
$$(X,Y)~\mapsto~\Lambda(X,Y)~=~\bigl(t(X,Y),~x(X,Y)\bigr).$$
\v
{\bf 6.} To complete the proof, 
we study in more detail the images of the singular sets
$S^w,S^z$.

By (\ref{impl1}) there can be only finitely many 
points inside 
$\Gamma$
where $w=\pi$ and $w_X=0$, say $P_i=(X_i,Y_i)$, $i=1,\ldots,m$. 
Moreover, by (\ref{impl2}), at a point $(X_0,Y_0)\in S^w\cap S^z$
we have
$$w_X\not=0, \qquad w_Y~=~0,\qquad z_X~=~0,\qquad z_Y~\not=~0.$$
Therefore, as shown in Fig.~\ref{f:wa67}, the two curves $\{w=\pi\}$
and $\{z=\pi\}$  intersect 
perpendicularly.    As a consequence, inside the compact set $\Gamma$, there can be only finitely many
such intersection points, say $Q_i = (X_i', Y_i')$, $i=1,\ldots,n$.

After removing these finitely many points $P_i, Q_i$,
we can thus write $\S^w$ as a finite union of curves $\gamma_j$
of the form
\bel{gj}\gamma_j~=~\bigl\{(X,Y)\,;~~X=\phi_j(Y),~~a_i<Y<b_j\bigr\}.\eeq
for suitable functions $\gamma_j$ of class $\C^2$.
We claim that the image of $\Lambda(\gamma_j)$ is a $\C^2$ curve
in the $t$-$x$ plane.  To prove this, it suffices to show that,
on the open interval $]a_j, b_j[\,$ the differential of the map
$$Y~\mapsto~\bigl( x(\phi_j(Y),Y),~t(\phi_j(Y),Y)\bigr)$$
does not vanish.
This is true because, by (\ref{4.2})
$${d\over dY} \, t(\phi_j(Y),Y)~=~t_X\cdot \phi_j' + t_Y
~=~0\cdot \phi_j'  + {(1+\cos z)q\over 4c(u)}~>~0.$$
Indeed, $z\not= \pi$ while $c(u), q >0$.

As shown in Fig.~\ref{f:wa68}, restricted to the 
inner domain $[0,T]\times [-\rho,\,\rho]$ in the $t$-$x$ plane, the singular set $\Lambda(S^w)$
is thus the union of the finitely many points 
$$\bega{l} p_i~=~\Lambda(P_i), \qquad i=1,\ldots,m,\\[4mm]
q_i~=~\Lambda(Q_i), \qquad i=1,\ldots,n,\enda$$
together with  finitely many $\C^2$ curves
$\Gamma_j~\doteq~\Lambda(\gamma_j)$. 
The same representation is valid for the image $\Lambda(S^z)$.
This concludes the proof of Theorem~1. \endproof
\v
\section{One-parameter families of solutions}
\setcounter{equation}{0}
 
In this section we study families of conservative solutions
$u=u(t,x,\lambda)$ of (\ref{1.1}) depending on an additional parameter
$\lambda\in [0,1]$.  
We thus consider a 1-parameter family of initial data
\bel{6.1}
u(0,x,\lambda)~=~u_0(x,\lambda),\qquad\qquad u_t(0,x,\lambda)~=~u_1(x,\lambda),\eeq
smoothly depending on the additional parameter $\lambda\in [0,1]$.
More precisely, these paths of initial data will lie in the space
\bel{TD}{\mathcal X}~\doteq~ 
\Big(\C^3([0,1]\times \R)\cap\L^\infty([0,1]; H^1(\R))\Big) 
\times\Big(\C^2([0,1]\times\R)\cap\L^\infty([0,1];\,\L^2(\R))\Big).\eeq
In particular, the map $(x,\lambda)\mapsto u_0(x,\lambda)$ is 
three times continuously differentiable and the $H^1$ norm
of $u_0(\cdot,\lambda)$ is uniformly bounded for all $\lambda$.
Moreover, the map $(x,\lambda)\mapsto u_1(x,\lambda)$ is 
two times continuously differentiable and the $\L^2$ norm
of $u_1(\cdot,\lambda)$ is uniformly bounded for all $\lambda$.

By an adaptation of the previous arguments one obtains
\v
{\bf Theorem 2.} 
{\it Let the wave speed $c(u)$ satisfy the assumptions {\bf (A)}
and let $T>0$ be given.
Then, for any 1-parameter family of initial data 
$(\hat u_0,\hat u_1)\in{\mathcal X}$
and any $\ve>0$, there exists a perturbed family
$(x,\lambda)\mapsto( u_0, u_1)(x,\lambda)$
such that
\bel{ed}
\Big\|( u_0- \hat u_0\,,~ u_1 - \hat u_1)\Big\|_{\mathcal X}
~<~\ve\,,\eeq
and moreover the following holds.
For all  except at most finitely many $\lambda\in [0,1]$, 
the conservative solution 
$u=u(t,x;\lambda)$ of (\ref{1.1})
is smooth in the complement of finitely
many points $P_i$ and finitely many $\C^2$ curves 
$\gamma_j$ in the domain $[0,T]\times \R$.}
\v
Toward a proof, we shall need
\v
{\bf Lemma 6.} {\it Let the function $u\mapsto c(u)$ satisfy the assumptions
{\bf (A)}, and let any $M>0$ be given.

Then there exists a dense set of paths of initial data 
$ \D~\subset~{\mathcal X}$ 
such that, if $(x,\lambda)\mapsto( u_0, u_1)(x,\lambda)$
lies in $\D$, then the corresponding solutions 
$(t,x,u,w,z,p,q)$ of (\ref{2.24})--(\ref{4.2})
with boundary data as in (\ref{2.28}) have the following properties.
On the domain $\Gamma$ in (\ref{GTdef})
one has
\begi\item[(i)] The map $(X,Y,\lambda)\mapsto (w,w_X,w_{XX})$
is transversal to the point $(\pi,0,0)$.
\item[(ii)] The map $(X,Y,\lambda)\mapsto (z,z_Y,z_{YY}) $
is transversal to the point $(\pi,0,0)$.
\item[(iii)] The map $(X,Y,\lambda)\mapsto (w,z,w_X)$
is transversal to the point $(\pi,\pi,0)$.
\item[(iv)] The map $(X,Y,\lambda)\mapsto (w,z,z_Y)$
is transversal to the point $(\pi,\pi,0)$.
\item[(v)] The map $(X,Y,\lambda)\mapsto (w,w_X,c'(u))$
is transversal to the point $(\pi,0,0)$.
\item[(vi)] The map $(X,Y,\lambda)\mapsto (z,z_Y,c'(u))$ 
is transversal to the point $(\pi,0,0)$.
\endi}
\v
{\bf Proof of Lemma~6.} Consider any point $(X_0, Y_0, \lambda_0)$.
Then,  there exist
3-parameter families of perturbed initial data
$(u_0^\theta, u_1^\theta)$, $\theta\in \R^3$ such that
the properties {\bf (1)--(3)} in Lemma 4 hold.
Indeed, it suffices to repeat all the arguments in the proof of Lemma~4
regarding $\lambda_0$ as a constant.   For a fixed 
$\lambda=\lambda_0$, the perturbations 
in (\ref{p00}) are thus functions of $s$ only, constant 
w.r.t.~$\lambda$.  

Combining these perturbations, as in the proof of 
Lemma~5, we obtain a map $(X,Y,\lambda, \theta)\mapsto (u,w,z,p,q)$
for which all transversality conditions (i)--(vi) are satisfied.
By the transversality theorem, for a.e.~$\theta$ the corresponding
map $(X,Y,\lambda)\mapsto 
(u^\theta,w^\theta,z^\theta,p^\theta,q^\theta)$ satisfies the
same transversality conditions.  This achieves the proof. \endproof

{\bf Proof of Theorem 2.}  
As in the proof of Theorem 1, we first 
choose $\rho$ large enough so that all our solutions
will be $\C^2$ for $(t,x)$ in the outer domain 
(\ref{outer}).  
 
For each $\lambda\in [0,1]$, we denote by
$(u,w,z,p,q, x,t)(X,Y,\lambda)$ the corresponding solution
of the semilinear system (\ref{2.26})-(\ref{4.2}).
We choose $M$ sufficiently large such that,
for all $\lambda\in [0,1]$, the inner domain 
$[0,T]\times [-\rho,\,\rho]$ is contained  in the image 
$$\Lambda^\lambda(\Gamma)~=~\Big\{\bigl( t(X,Y,\lambda)
,~x(X,Y,\lambda)\bigr)\,;~~|X|+|Y|\leq M\Big\}.$$

By performing an arbitrarily 
small perturbation of the initial path of solutions
we obtain a second path $\lambda\mapsto u(\cdot,\lambda)$
such that, in  the corresponding solution 
$(u,w,z,p,q, x,t)(X,Y,\lambda)$,
the transversality relations (i)--(vi) in Lemma~6 hold.

\begin{figure}[htbp]
\centering
\includegraphics[width=1.0\textwidth]{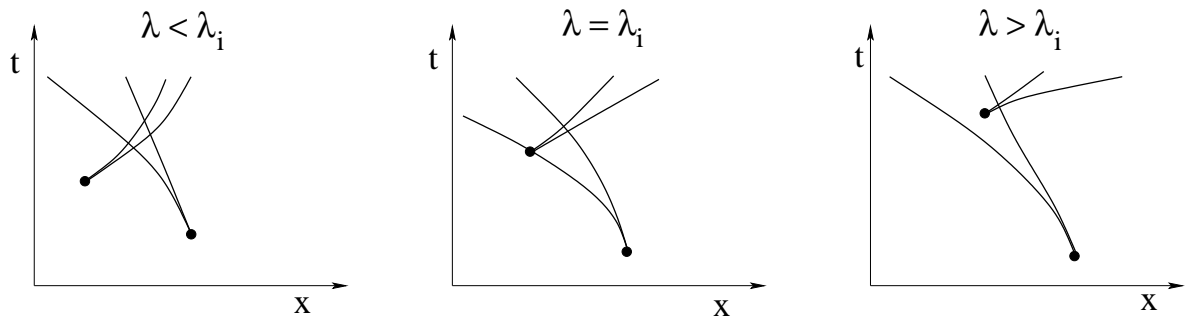}
\caption{ \small The singular set for a solution $u(t,x;\lambda)$.
When the parameter $\lambda$ crosses one of the critical 
values $\lambda_i$, the topology of the singular set  changes.}
\label{f:wa70}
\end{figure}

Since the variables
$(X,Y,\lambda)\in \Gamma\times [0,1]$ range on a compact,
three dimensional set,  this implies that the 
values in (i)--(vi) are attained only at finitely many points,
say $(X_i, Y_i,\lambda_i)$, $i=1,\ldots,n$.
Hence, for $\lambda\notin \{\lambda_1,\ldots,\lambda_n\}$, 
 the solution $(t,x,u,w,z,p,q)(\cdot,\cdot,\lambda)$ 
 does not attain any of the  values in (i)--(vi),
 for $(X,Y)\in \Gamma$.   
As shown in steps
{\bf 5-6} in the proof of Theorem 1, 
the  corresponding solution 
 $u=u(t,x;\,\lambda)$
is then piecewise smooth on the inner domain 
$[0,T]\times [-\rho,\,\rho]$.  \endproof
\v
{\bf Remark 5.} For a given solution $u=u(t,x)$, 
define its {\it singular set} as
$$S^u~\doteq~\bigl\{ (t,x)\,;~~ u~\hbox{is not $\C^2$ on any
 neighborhood of}~ (t,x)\bigr\}.$$
In the above construction, one 
can regard $\lambda_1,\ldots,\lambda_n$
as {\it bifurcation values}, where the structure of the singular set
changes (Fig.~\ref{f:wa70}). 
On the other hand, for $\lambda\notin\{\lambda_1,\ldots,\lambda_n\}$ 
the solution $u(\cdot,\cdot;\lambda)$ is 
{\it structurally stable}.   A small perturbation of the initial data
does not change the topology of the singular set.
Based on the present analysis, we speculate that a theory of 
generic structural stability and a global classification
of solutions to (\ref{1.1}) can be developed, in analogy to 
the classical theory for ODEs \cite{A, P}.

\v
{\bf Acknowledgments.} This research was partially supported by NSF, with grant DMS-1411786: Hyperbolic Conservation Laws and Applications, and by the AMS Simons Travel Grant Program.
\v

\end{document}